\documentclass[11pt]{amsart}
\usepackage{amssymb,amsmath,amsthm,amscd}

\input{amssym.def}

\newcommand{\bt}{\begin{Theorem}}
\newcommand{\et}{\end{Theorem}}
\newcommand{\bi}{\begin{itemize}}
\newcommand{\ei}{\end{itemize}}
\newcommand{\bea}{\begin{eqnarray}}
\newcommand{\eea}{\end{eqnarray}}

\theoremstyle{plain}

\newtheorem{Theorem}{\sc Theorem}[section]

\newtheorem{Proposition}{\sc Proposition}[section]

\theoremstyle{definition}

\theoremstyle{remark}

\newcommand{\be}{\begin{equation}}
\newcommand{\ee}{\end{equation}}

\newcommand{\inner}[2]{\langle #1,#2 \rangle }%
\title[classification of homogeneous operators]{A  classification of homogeneous operators\\ in the Cowen-Douglas class}

\author{Adam Kor\'{a}nyi}
\address{Lehman College\\
The City University of New York\\
Bronx, NY 10468
}
\email{adam.koranyi@lehman.cuny.edu}
\author{Gadadhar Misra}
\address{Department of Mathematics \\
Indian Institute of Science\\
Bangalore 560 012
 }
\email{gm@math.iisc.ernet.in}
\thanks{This research was supported in part by a DST -
NSF S\&T Cooperation Programme and a PSC-CUNY grant.}
\begin{document}
\maketitle
\begin{abstract}
A complete list of homogeneous operators in the Cowen - Douglas class ${\rm B}_n(\mathbb D)$ is given.  This classification is obtained from an explicit realization of all the homogeneous Hermitian  holomorphic vector bundles on the unit disc under the action of the universal covering group of the bi-holomorphic automorphism group of the unit disc.
\end{abstract}
\section{Introduction}
\noindent
A bounded linear  operator $T$ on a complex separable Hilbert space $\mathcal H$ is said to be {\em homogeneous} if its spectrum is contained in the closed unit disc and for every M\"{o}bius transformation $g$ of the unit disc $\mathbb D$, the operator $g(T)$ defined via the usual holomorphic functional calculus, is unitarily equivalent to $T$.  To every homogeneous irreducible operator $T$ there corresponds an {\em associated projective unitary representation}  $U$ of the  M\"{o}bius group $G$:
$$U_g^*\, T \, U_g =  g(T),\: g \in G.$$
The projective unitary representations of $G$ lift to unitary representations of the universal cover $\tilde{G}$ which are quite well-known.  We can choose $U_g$ such that $ k\mapsto U_k$ is a representation of the rotation group $\mathbb K \subseteq G$. 
If
$$
\mathcal H(n) = \{ x  \in \mathcal H: U_k x = e^{i\,n \theta} x\},
$$
then $T:\mathcal H(n) \to \mathcal H(n+1)$ is a block shift.
A complete classification of these for $\dim \mathcal H(n) \leq 1$ was obtained in \cite{shift} using the representation theory of $\tilde{G}$.  First examples for $\dim \mathcal H(n) = 2$ appeared in \cite{W}.  Recently (cf. \cite{KM, KMISI}), a $m$ - parameter family of examples with $\dim \mathcal H(n) = m$ was constructed.
In the present announcement we show that the ideas of \cite{KM, KMISI} lead to a complete classification of the homogeneous operators in the Cowen - Douglas class.

A Fredholm operator $T$ on a Hilbert space $\mathcal H$ is said to be in the Cowen - Douglas class of the domain $\Omega \subseteq \mathbb C$ if its eigenspaces $E_w,\, w\in \Omega$ are of constant finite dimension.  In the paper \cite{C-D}, Cowen and Douglas show that
\begin{enumerate}
\item[(a)] $E \subseteq \Omega \times \mathcal H$ with fiber $E_w$ at $w\in \Omega$ is a holomorphic Hermitian vector bundle over $\Omega$, where the Hermitian structure is given by
$$
\|s_w\|_w = \|\iota_w s_w\|_\mathcal H,\,s_w \in E_w,
$$
and $\iota_w: E_w \to \mathcal H$ is the inclusion map;
\vskip 0.5ex
\item[(b)] isomorphism classes of $E$ correspond to unitary equivalence classes of $T$;\vskip 0.5ex
\item[(c)] the holomorphic Hermitian vector bundle $E$ is irreducible if and only if the operator $T$ is irreducible.
\end{enumerate}

It can be shown that a Cowen-Douglas operator is homogeneous if and only if the corresponding bundle is homogeneous under $\tilde{G}$.  We describe below all irreducible homogeneous holomorphic Hermitian vector bundles over the unit disc and determine which ones of these correspond to homogeneous operators (necessarily irreducible) in the Cowen-Douglas class.

\section{Homogeneous holomorphic vector bundles}
The description of homogeneous vector bundles via holomorphic induction is well-known.
Let $\mathfrak t \subseteq \mathfrak g^\mathbb C=\mathfrak s\mathfrak l(2,\mathbb C)$ be the algebra $ \mathbb C h + \mathbb C y$, where
$$ h = \frac{1}{2}\Big (\begin{matrix} 1 & 0\\ 0 & - 1 \end{matrix} \Big ),\, y = \Big (\begin{matrix} 0 & 0\\ 1 & 0 \end{matrix} \Big ).$$
Linear representations $(\varrho, V)$ of the algebra $\mathfrak t \subseteq \mathfrak g^\mathbb C=\mathfrak s\mathfrak l(2,\mathbb C)$, that is, pairs $\varrho(h), \varrho(y)$ of linear transformations satisfying $[\varrho(h),\varrho(y)] = - \varrho(y)$ provide a para-metrization of the homogeneous holomorphic vector bundles.

The $\tilde{G}$ - invariant Hermitian structures on the homogeneous holomorphic vector bundle $E$ (making it into a homogeneous holomorphic Hermitian vector bundle), if they exist, are given by $\varrho(\tilde{\mathbb K})$ - invariant inner products on the representation space.  Here $\tilde{\mathbb K}$ is the stabilizer of $0$ in $\tilde{G}$.

An inner product can  be $\varrho(\tilde{\mathbb K})$ - invariant if and only if $ \varrho(h)$ is diagonal with real diagonal elements in an appropriate basis.  We are interested only in Hermitizable bundles, that is, those that admit a Hermitian structure.  So, we will assume without restricting generality, that the representation space of $\varrho$ is $\mathbb C^n$ and that $\varrho(h)$ is a real diagonal matrix.

Since $[\varrho(h),\varrho(y)] = - \varrho(y)$, we have $\varrho(y) V_\lambda \subseteq V_{\lambda -1}$, where $V_\lambda = \{ \xi \in \mathbb C^n: \varrho(h) \xi = \lambda \xi\}$.
Hence $(\varrho, \mathbb C^n)$ is a  direct sum, orthogonal for every $\varrho(\tilde{K})$ - invariant inner   product of ``elementary'' representations,  that is, such that
$$ \varrho(h) = \begin{pmatrix}
-\eta I_0 & &\\
 &\ddots & \\
&&-(\eta + m) I_m
\end{pmatrix}\mbox{with}\,\, I_j = I\,\, \mbox{on}\,\, V_{-(\eta+j)} = \mathbb C^{d_j}$$
and
$$ Y:= \varrho(y) = \begin{pmatrix}
0 & & & \\
Y_1 & 0 & &\\
&Y_2&0 & &\\
&&\ddots&\ddots&\\
&&&Y_m&0
\end{pmatrix},\, Y_j:V_{-(\eta+j-1)} \to V_{-(\eta+j)}.$$
We denote the corresponding elementary Hermitizable bundle by $E^{(\eta, Y)}$. 

\subsection{\sc The Multiplier and Hermitian  structures}
As in \cite{KMISI} we will use a natural trivialization of $E^{(\eta, Y)}$. In  this  the sections of homogeneous holomorphic vector bundle $E^{(\eta, Y)}$ are holomorphic functions $\mathbb D$ taking values in $\mathbb C^n$.  The $\tilde{G}$ action is given by $f  \mapsto J_{g^{-1}}^{(\eta,  Y)} \big (f\circ g^{-1}\big )$ with multiplier
$$(\!\!( J^{(\eta, Y)}_g(z))\!\!)_{p,\ell} =\begin{cases}
\tfrac{1}{(p-\ell)!} (-c_g)^{p-\ell}
(g^\prime)(z)^{\eta + \frac{p + \ell}{2}} Y_p \cdots Y_{\ell+1} &{~if~} p \geq \ell\\
0&{~if~} p < \ell \end{cases},
$$
where $c_g$ is  the analytic function on $\tilde{G}$ which, for $g$ near $e$, acting on $\mathbb D$ by $z \mapsto (a z + b) (cz +d)^{-1}$ agrees with $c$.
\begin{Proposition}
We have $E^{(\eta, Y)}\equiv E^{(\eta^\prime, Y^\prime)}$ if and only if $\eta=\eta^\prime$ and $Y^\prime = AYA^{-1}$ with a block diagonal matrix $A$.
\end{Proposition}
A Hermitian structure on $E^{(\eta, Y)}$ appears as the assignment of an inner product  $\langle \cdot , \cdot \rangle_z$ on  $\mathbb C^n$ for   $z \in \mathbb D$. We can write 
$$ \langle \zeta, \xi \rangle_z = \langle H(z) \zeta , \xi \rangle,\,\,\mbox{with}\,\, H(z) \succ 0.$$
Homogeneity as a Hermitian vector bundle is equivalent to
$$ J_g(z) H(g\cdot z)^{-1}J_g(z)^* = H(z)^{-1},\,\,g \in G,\,\, z\in \mathbb D.$$
The Hermitian   structure is then determined by $H=H(0)$ which is a positive block diagonal matrix.  We write $(E^{(\eta, Y)}, H)$ for the vector bundle $E^{(\eta, Y)}$ equipped with the
Hermitian structure $H$.  We note that $(E^{(\eta, Y)}, H) \cong (E^{(\eta, AYA^{-1})}, {A^*}^{-1}HA)$ for any block diagonal invertible $A$.  Therefore every homogeneous holomorphic Hermitian vector bundle is isomorphic with one of the form  $(E^{(\eta, Y)}, I)$.

If $E^{(\eta, Y)}$ has a reproducing kernel $K$ which
is the case for bundles corresponding to an operator in the  Cowen-Douglas class, 
then $K$ satisfies 
$$K(z,w) = J_g(z) K(g z , g w) J_g(w)^*$$
and induces a Hermitian structure $H$ given by $H(0) = K(0,0)^{-1}.$
\section{Construction of the bundles with reproducing kernel}
For $\lambda >0$, let $\mathbb A^{(\lambda)}$ be the Hilbert space of holomorphic functions on the unit disc with reproducing kernel $(1-z \bar{w})^{-2 \lambda}$.
It is  homogeneous under the multiplier  $g^\lambda(z)$ for the action of $\tilde{G}$. This gives a unitary representation of $\tilde{G}$.  Let  $\mathbf A^{(\eta)}= \oplus_{j=0}^m \mathbb A^{(\eta + j)} \otimes \mathbb C^{d_j}$.  For $f$ in $\mathbf A^{(\eta)}$, we denote by $f_j$, the part of $f$ in $\mathbb A^{(\eta + j)} \otimes \mathbb C^{d_j}$.  We define $\Gamma^{(\eta, Y)} f$ as the $\mathbb C^n$ - valued holomorphic function whose part in $\mathbb C^{d_\ell}$ is given by
$$ \big (\Gamma^{(\eta, Y)} f \big )_\ell = \sum_{j=0}^\ell\frac{1}{(\ell -j)!}\frac{1}{(2 \eta + 2 j)_{\ell -j}} Y_\ell \cdots Y_{j+1} f_j^{(\ell -j)}$$
for $\ell \geq j$. For invertible block diagonal $N$ on $\mathbb C^n$, we also define 
$\Gamma_N^{(\eta, Y)}:= \Gamma^{(\eta, Y)} \circ N$.  It can be verified that $ \Gamma_N^{(\eta, Y)}$ is a $\tilde{G}$ -   equivariant isomorphism of $\mathbf A^{(\eta)}$ as a homogeneous holomorphic vector bundle onto $E^{(\eta, Y)}$.  The image $K_N^{(\eta, Y)}$ of the reproducing kernel of $\mathbf A^{(\eta)}$ is then a reproducing kernel for $E^{(\eta, Y)}$.  A computation gives that $K^{(\eta,Y)}_N(0,0)$ is a block diagonal matrix such that its $\ell$'th block is
$$
K^{(\eta, Y)}_N(0,0)_{\ell,\ell} = \sum_{j=0}^\ell \frac{1}{(\ell -j)!}\frac{1}{(2\eta + 2j)_{\ell-j}} Y_\ell\cdots Y_{j+1} N_jN_j^* Y_{j+1}^* \cdots Y_\ell^*.
$$
We set $H^{(\eta,Y)}_N = {K_N^{(\eta, Y)}(0,0)}^{-1}$.  We have now constructed a family $(E^{(\eta, Y)}, H^{(\eta,Y)}_N)$ of elementary homogeneous holomorphic vector bundles with a reproducing kernel  ($\eta >0$, $Y$ as before, $N$ invertible block diagonal).
\begin{Theorem}
Every elementary homogeneous holomorphic vector bundle $E$ with a reproducing kernel arises
from the construction given above.
\end{Theorem}
\begin{proof}[Sketch of proof]
As a homogeneous bundle $E$ is isomorphic to some $E^{(\eta, Y)}$.   Its reproducing kernel gives a Hilbert space structure in which the $\tilde{G}$ action on the sections of $E^{(\eta, Y)}$ is a unitary representation $U$.  Now $\Gamma^{(\eta, Y)}$ intertwines the unitary representation of $\tilde{G}$ on $\mathbb A^{(\eta)}$ with $U$.  The existence of a block diagonal $N$ such that $\Gamma_N^{(\eta, Y)} = \Gamma^{(\eta, Y)} \circ N$ is a Hilbert space isometry follows from Schur's Lemma.
\end{proof}

As remarked before, every homogeneous holomorphic Hermitian vector bundle is isomorphic to an $(E^{(\eta,Y)}, I)$, here $Y$ is unique up to conjugation by a block unitary.  In this form, it is easy to tell whether the bundle is irreducible:  this is the case if and only if $Y$ is not the orthogonal direct sum of two matrices of the same block type as $Y$. We call such a $Y$ irreducible.

Let $\mathcal P$ be the set of all $(\eta, Y)$ such that $E^{(\eta, Y)}$ has a reproducing kernel.  Using the formula for $K_N^{(\eta, Y)}(0,0)$ we can write down explicit systems of inequalities that determine whether $(\eta, Y)$ is in $\mathcal P$.    In particular we have
\begin{Proposition}
For every $Y$, there exists a $\eta_Y > 0$ such that $(\eta, Y)$ is in $\mathcal P$ if and only if $\eta>\eta_Y $.
\end{Proposition}

Finally, we obtain the announced classification.
\begin{Theorem}
All the homogeneous holomorphic Hermitian vector bundles of rank $n$ with a reproducing kernel correspond to homogeneous operators in the Cowen -- Douglas class ${\rm B}_n(\mathbb D)$. The irreducible ones are the adjoint of the multiplication operator $M$ on  the Hilbert space of sections of $(E^{(\eta, Y)},I)$ for some $(\eta, Y)$ in $\mathcal P$ and irreducible $Y$.  The block matrix $Y$ is determined up to conjugacy by block diagonal unitaries.
\end{Theorem}

\begin{proof}[Sketch of proof]
There is a simple orthonormal system for the Hilbert space  $\mathbb A^{(\lambda)}$.  Hence we can find such a system for $\mathbf A^{(\eta)}$ as well.  Transplant it using  $\Gamma^{(\eta, Y)}$ to   $E^{(\eta,  Y)}$.
The multiplication operator in this basis has a block diagonal form with  { $M_n:= M_{|\,{\rm res\,\,} \mathcal H(n)}: \mathcal H(n) \to \mathcal H(n+1)$.} This description is sufficiently explicit to see:
{$ M_n \sim I + 0(\frac{1}{n})$.} Hence {$M$} is the sum of an ordinary  block shift operator and a Hilbert Schmidt operator. This completes the proof.
\end{proof}

\bibliographystyle{amsplain}
\providecommand{\bysame}{\leavevmode\hbox to3em{\hrulefill}\thinspace}
\providecommand{\MR}{\relax\ifhmode\unskip\space\fi MR }
\providecommand{\MRhref}[2]{%
  \href{http://www.ams.org/mathscinet-getitem?mr=#1}{#2}
}
\providecommand{\href}[2]{#2}

\end{document}